\documentclass[11pt]{article}
\usepackage{amsfonts}
\usepackage{amsfonts}
\usepackage{amsfonts}
\usepackage{mathrsfs}
\usepackage{amssymb,amsmath}
 \allowdisplaybreaks

\catcode`!=11
\let\!int\int \def\int{\displaystyle\!int}
\let\!lim\lim \def\lim{\displaystyle\!lim}
\let\!sum\sum \def\sum{\displaystyle\!sum}
\let\!sup\sup \def\sup{\displaystyle\!sup}
\let\!inf\inf \def\inf{\displaystyle\!inf}
\let\!cap\cap \def\cap{\displaystyle\!cap}
\let\!max\max \def\max{\displaystyle\!max}
\let\!min\min \def\min{\displaystyle\!min}
\let\!frac\frac \def\frac{\displaystyle\!frac}
\catcode`!=12

\let\oldsection\section
\renewcommand\section{\setcounter{equation}{0}\oldsection}

\allowdisplaybreaks
\def\pf{\it{Proof.}\rm\quad}

\newcommand\dd{\mbox{d}}

\def\R{\mathbb{R}}

\newcommand\divg{{\text{div}}}

\newtheorem{thm}{Theorem}[section]
\newtheorem{pro}{Proposition}[section]
\newtheorem{lem}{Lemma}[section]

\begin{document}
\title{\bf Global well-posedness for the
two-dimensional equations of nonhomogeneous incompressible liquid crystal flows
with nonnegative density\thanks{This work was partially supported by
  NNSFC (Grant Nos.  11271306 \& 10971171),
  the Fundamental Research Funds for the Central Universities (Grant Nos. 2010121006 \& 2012121005),
  and the Natural Science Foundation of Fujian Province of China (Grant No.
  2010J05011).
      }}
\author{ Shengquan Liu$^1$, Jianwen Zhang$^{2,}$\thanks{Corresponding author (E-mail:jwzhang@xmu.edu.cn).}\\[2mm]
\small $^1$ School of Mathematics, Liaoning University, Shenyang 110036, China, \\
\small $^2$ School of Mathematical Sciences, Xiamen University,
Xiamen 361005, China}
\date{}
\maketitle \noindent{\bf Abstract.} In this paper, the authors first
establish the global well-posedness of strong solutions of the
simplified Ericksen-Leslie model for nonhomogeneous incompressible
nematic liquid crystal flows in two dimensions if the initial data
satisfies some smallness condition. It is worth pointing out that
the initial density is allowed to contain vacuum states and the
initial velocity can be arbitrarily large. We also present a Serrin's type criterion,
depending only on $\nabla d$, for the breakdown of local strong
solutions. As a byproduct, the global strong solutions with large initial data are
obtained, provided the macroscopic molecular orientation of the liquid crystal materials
satisfies a natural geometric angle condition (cf. \cite{LLZ2012}).

\vskip 2mm

\noindent{\bf Keywords.} Liquid crystals; Nonhomogeneous incompressible flows; Global strong solutions; Vacuum; Blowup criterion
\\[2mm]
\noindent{\bf AMS Subject Classifications (2000).} 35B45, 76A15, 76D03, 76D05


\section{Introduction}
Liquid crystals are substances that exhibit a phase of matter that has properties between those of a conventional liquid and
those of a solid crystal (cf. \cite{Ge1974}).
The hydrodynamic theory of liquid crystals was first developed by Ericken and Leslie during the period of 1958 through 1968
(see \cite{Er1961,Er1962,Le1968,Le1979}).
Since then, many remarkable developments have been made from both theoretical and applied aspects, however, many physically
important and mathematically fundamental problems still remain open.  In this paper, we
consider a simplified Ericken-Leslie model for the nonhomogeneous incompressible nematic liquid crystals in two dimensions:
\begin{eqnarray}
\rho_t+\divg(\rho u)&=&0,  \label{1.1}\\
(\rho u)_t+\divg(\rho u\otimes u)+\nabla P&=&\Delta u- \nabla d\cdot \Delta d,\label{1.2}\\
\divg
u&=&0,\label{1.3}\\
d_t+ u\cdot\nabla d&=&\Delta d+|\nabla d|^2d,\label{1.4}
\end{eqnarray}
where $\rho : \R^2\times[0,\infty)\to \R^+$ is the density of the
fluid, $u :  \R^2\times[0,\infty)\to\R^2$ is the velocity field of
the fluid, $P :  \R^2\times[0,\infty)\to\R $ is the pressure of the
fluid, and $d :  \R^2\times[0,\infty)\to \mathbb{S}^2$ (the unit
sphere in $\R^3$, i.e. $|d|=1$) represents the averaged
macroscopic/continuum molecular orientations.

Though system (\ref{1.1})--(\ref{1.4}) is a simplified version of
the Ericksen-Leslie model, but it still retains the most
interesting mathematical properties without losing the basic
nonlinear structure of the original Ericksen--Leslie model
\cite{Er1961,Er1962,Le1968,Le1979}. Roughly speaking, the system
\eqref{1.1}--\eqref{1.4} is a system of the nonhomogeneous
Navier-Stokes equations for incompressible flows coupled with the
equation for heat flow of harmonic maps, and thus, its mathematical
analysis is full of challenges. In particular, if $\rho={\rm
Const.}$, then it turns into the following homogeneous system which
models the incompressible flows of nematic liquid crystal
\begin{eqnarray}
u_t+u\cdot\nabla u+\nabla P&=&\mu\Delta u- \nabla d\cdot \Delta d,\label{1.5}\\
\divg u&=&0,\label{1.6}\\
d_t+ u\cdot\nabla d&=&\Delta d+|\nabla d|^2d\label{1.7}
\end{eqnarray}
with $|d|=1$.
Moreover, if $u=0$ in (\ref{1.5})--(\ref{1.7}), then it reduces to the following equation for heat flow of harmonic maps:
\begin{equation}
d_t=\Delta d+|\nabla d|^2d,\quad |d|=1.\label{1.8}
\end{equation}

There has been a lot of literature on the mathematical studies of
(\ref{1.5})--(\ref{1.7}) and (\ref{1.8}), see, for example,
\cite{HW2012,Hong2011,Lin1989,LL1995,LL1996,LLW2010,XZ2012,Wang2011}
and \cite{Chang1989,CD1990,CDY1992,CS1989,St1985}, respectively. In
the following, we briefly recall some related mathematical results
of the liquid crystal flows. In a series of papers, Lin
\cite{Lin1989} and Lin-Liu \cite{LL1995,LL1996} initiated the
mathematical analysis of (\ref{1.5})--(\ref{1.7}) in 1990s. More
precisely, to relax the nolinear constraint $|d|=1$, they proposed
an approximate model of Ericksen-Leslie system with variable length
by Ginzburg-Landau functionals, that is, the equation (\ref{1.7})
with $|d|=1$ is replaced by
\begin{equation}
d_t+u\cdot\nabla d=\Delta d+\frac{1}{\varepsilon^2}\left(1-|d|^2\right)d.\label{1.9}
\end{equation}
In \cite{Lin1989,LL1995}, the authors proved the global existence of
classical and weak solutions of (\ref{1.5}), (\ref{1.6}),
(\ref{1.9}) in dimensions two and three, respectively. The partial
regularity of suitable weak solutions was also studied in
\cite{LL1996}. However, as pointed out in \cite{LL1995}, the
vanishing limit of $\varepsilon\to0$ is an open and challenging
problem. Indeed, in contrast with  (\ref{1.9}), it is much more
difficult to deal with the nonlinear term $|\nabla d|^2d$ with
$|d|=1$ appearing on the right-hand side of (\ref{1.4}) or
(\ref{1.7}) from the mathematical point of view. In two independent
papers \cite{Hong2011} and \cite{LLW2010}, Hong and Lin-Lin-Wang
showed the global existence of weak solutions of
(\ref{1.5})--(\ref{1.7}) in dimensions two, and proved that the
solutions are smooth away from at most finitely many singular times
which is analogous to that for the heat flows of harmonic maps (see
\cite{Chang1989,St1985}). The global existence of smooth solution
with small initial data of (\ref{1.5})--(\ref{1.7}) was also proved
\cite{LLW2010,XZ2012} and \cite{Wang2011,LW2012} in dimensions two and
three, respectively.

For the approximate nonhomogeneous equations
(\ref{1.1})--(\ref{1.3}) and (\ref{1.9}), the global existence of
weak solutions with generally large initial data was proved in
\cite{LZ2009,JT2009}, and the global regularity of the solution with
strictly positive density was studied in \cite{Dai2012}. As
aforementioned, the nonlinear term $|\nabla d|^2d$ with $|d|=1$ will
cause serious difficulty in the mathematical analysis of liquid
crystal flows. Recently, Wen and Ding \cite{WD2011} established the
local existence and uniqueness of strong solutions of
(\ref{1.1})--(\ref{1.4}) in the case that the initial density may
contain vacuum states (i.e. $\rho_0\geq0$). Moreover, if the initial
density has a positive lower bound (i.e. $\rho_0\geq
\underline\rho>0$) which indicates that there is absent of vacuum
initially, the global strong solutions with small initial data was
also obtained in \cite{WD2011}.

As that for the density-dependent Navier-Stokes equations (see
\cite{CK2003,Li1996}), the possible presence of vacuum is one of the
major difficulties when the problems of global existence, uniqueness
and regularity of solutions are involved. Therefore, in the present paper we aim to
investigate the global regularity of (\ref{1.1})--(\ref{1.4}) when the initial density may contain vacuum.

We consider the Cauchy problem of
\eqref{1.1}--\eqref{1.4} with the following initial data:
\begin{equation}
(\rho,u,d)(x,0)=(\rho_0,u_0,d_0)(x)\quad{\rm for}\quad  x\in\R^2,\label{1.10}
\end{equation}
and the far-field behavior at infinity:
\begin{equation}
\quad(\rho,u,d)(x,t)\to(\tilde\rho,0,e)\quad{\rm as}\quad |x|\to\infty,\;
t>0,\label{1.11}
\end{equation}
where $\tilde\rho>0$ is a given positive constant and $e\in \mathbb{S}^2$ is
a given unit vector (i.e. $|e|=1$).

To state our main results, we first introduce the definition of
strong solutions of (\ref{1.1})--(\ref{1.4}), (\ref{1.10}) and
(\ref{1.11}).

\vskip 2mm

\noindent{\bf Definition 1.1} A pair of functions $(\rho,u,P,d)$ is called a strong solution of (\ref{1.1})--(\ref{1.4}),
(\ref{1.10}) and (\ref{1.11})
on $\R^2\times[0,T]$, if $\rho(x,t)\geq0$ for all $(x,t)\in\R^2\times[0,T]$,
\begin{equation}\label{1.12}
\left\{
\begin{array}{lll}
&\rho-\tilde{\rho}\in C([0, T]; H^2(\R^2)),\quad \rho_t\in L^\infty(0,T;H^1(\R^2))\\[2mm]
&u\in C([0, T]; H^2(\R^2))\cap L^2(0,T; H^{3}(\R^2)),\\[2mm]
&\sqrt{\rho}u_t\in L^\infty(0, T; L^2(\R^2)),\quad u_t\in L^2(0, T; H^1(\R^2)),\\[2mm]
&\nabla P\in C([0, T]; L^2(\R^2))\cap L^2(0, T; H^1(\R^2)),\\[2mm]
&\nabla d\in C([0, T]; H^2(\R^2)),\quad d_t\in L^\infty(0, T; H^1(\R^2))\cap L^2(0, T; H^2(\R^2)),
\end{array}
\right.
\end{equation}
and $(\rho,u,P,d)$ satisfies (\ref{1.1})--(\ref{1.4}) a.e. on $\R^2\times(0,T]$.

\vskip 2mm

Then, our first result concerning the global strong solutions with
small data can be stated in the following theorem.
\begin{thm}\label{thm1.1}
Assume that the initial data $(\rho_0,u_0,d_0)$ satisfies
\begin{equation}
\label{1.13}
\left\{
\begin{array}{lll}
\rho_0\geq0,\quad (\rho_0-\tilde{\rho},u_0,\nabla d)\in
H^2(\R^2),\quad {\rm div}u_0=0,\quad |d_0|=1,\\[2mm]
\Delta u_0-\nabla P_0-\nabla d_0\cdot\Delta d_0=\rho_0^{1/2}g\quad{for\; some }\quad (\nabla P_0, g)\in L^2(\R^2).
\end{array}\right.
\end{equation}
Then for any given $0<T<\infty$, there exists a unique global strong solution $(\rho, u, P, d)$ of \eqref{1.1}--\eqref{1.4}, \eqref{1.10}
and \eqref{1.11} on $\R^2\times[0,T]$, provided
\begin{equation}
 \exp{\left(2\left(\|\rho_0^{1/2}u_0\|^2_{L^2}+\|\nabla d_0\|^2_{L^2}\right)\right)}\|\nabla d_0\|^2_{L^2}\leq \frac{1}{16}.\label{1.14}
\end{equation}
\end{thm}

It is worth mentioning that the smallness condition (\ref{1.14})
stated in Theorem \ref{thm1.1} implies that $(\rho_0,u_0)$ can be
arbitrarily large if $\|\nabla d_0\|_{L^2}$ is chosen to be
suitably small. This is analogous to the one in \cite{XZ2012}.
Moreover, as a result, we see that the strong solution to the Cauchy
problem of nonhomogeneous Navier-Stokes equations (i.e. $d={\rm
Const.}$) with large initial data, which may contain vacuum, exists
globally on $\R^2\times[0,T]$ for all $0<T<\infty$. Thus, Theorem
\ref{thm1.1} also generalizes the result due to Huang-Wang
\cite{HW2012}.

The proof of Theorem \ref{thm1.1} is mainly based on a critical
Sobolev inequality of logarithmic type which was recently proved by
Huang-Wang (cf. \cite{HW2013}) and is originally due to
Brezis-Wainger \cite{BW1980} (see also \cite{KOT2002,Oz1995}).
However, it is remarkable that the arguments in \cite{HW2013}
actually depend on the size of the domain considered and cannot be
applied directly to the case of the whole space. Thus, some new
ideas have to be developed. The main difference lies in the proof of
Lemma \ref{lem3.3}, where, instead of $\|\rho^{1/2} u_t\|_{L^2}$ and
$\|\rho^{1/2} u\cdot\nabla u\|_{L^2}$, we use the material
derivative $\|\rho^{1/2}\dot u\|_{L^2}$ for some technical reasons.
We also note here that the strictly positive far-field condition
$\tilde\rho>0$ plays an important role in our analysis. The strongly
nonlinear terms $|\nabla d|^2d$ and $\nabla d\cdot\Delta d$ in
(\ref{1.2}) and (\ref{1.4}) will also cause some additional
difficulties.

For the generally large initial data, it is still an interesting and
open problem whether the strong solution blows up or not in finite
time.  In \cite{LLW2010} and \cite{HW2012}, the authors proved
respectively that the following blowup criteria for the
two-dimensional equations of (\ref{1.5})--(\ref{1.7}):
\begin{equation}
\lim\limits_{T\to T^*}\int_0^T\left(\| u\|_{L^4}^4+\|\nabla
d\|_{L^4}^4\right)dt=\infty\quad{\rm and}\quad \lim\limits_{T\to
T^*}\int_0^T\|\nabla d\|_{L^\infty}dt=\infty,\label{1.15}
\end{equation}
where $0<T^*<\infty$ is the maximal time of the existence
of a strong solution to (\ref{1.5})--(\ref{1.7}).
Motivated by the proofs of Theorem \ref{thm1.1}, we can prove the
following  mechanism for possible breakdown of strong solutions, which is a natural extension of the ones in
\cite{LLW2010,HW2012}.

\begin{thm}\label{thm1.2}Assume that $0<T^*<\infty$ is the maximal time of the existence
of a strong solution to \eqref{1.1}--\eqref{1.4}, \eqref{1.10} and
\eqref{1.11} with generally large initial data $(\rho_0,u_0,d_0)$
satisfying \eqref{1.13}. Then,
\begin{equation}
\lim\limits_{T\to T^*}\int_0^T\|\nabla
d\|_{L^r}^sdt=\infty\label{1.16}
\end{equation}
for any $(r,s)$ satisfying
\begin{equation}
\frac{1}{r}+\frac{1}{s}\leq \frac{1}{2},\quad 2<
r\leq\infty.\label{1.17}
\end{equation}
\end{thm}

 Theorem \ref{1.2} implies that for any $0<T<\infty$ if the left-hand side of (\ref{1.16}) is finite, then the strong solution of  (\ref{1.1})--(\ref{1.4}), (\ref{1.10}) and
(\ref{1.11}) will exist globally on $\R^2\times(0,T)$.

Based on a frequency localization argument combined with the
concentration-compactness approach, Lei-Li-Zhang \cite{LLZ2012}
recently proved the following interesting rigidity theorem for the
approximate harmonic maps.
\begin{pro}(\cite[Theroem 1.5]{LLZ2012})\label{pro1.1} For given positive constants $0<C_0<\infty$ and $0<\varepsilon\leq 1$, assume that $d:\R^2\to \mathbb{S}^2$ satisfying $\nabla d\in H^1(\R^2)$ with $\|\nabla d\|_{L^2}\leq C_0$ and $d_3\geq\varepsilon$. Then there exists a positive constant $\delta_0\in(0,1)$, which depends only on $C_0$ and $\varepsilon$, such that
\begin{equation}
\|\nabla d\|_{L^4}^4\leq
\left(1-\delta_0\right)\|\nabla^2d\|_{L^2}^2,\label{1.18}
\end{equation}
which particularly implies
\begin{equation}
\|\Delta d+|\nabla
d|^2d\|_{L^2}^2\geq\frac{\delta_0}{2}\left(\|\Delta
d\|_{L^2}^2+\|\nabla d\|_{L^4}^4\right).\label{1.19}
\end{equation}
\end{pro}

As an immediate consequence of Theorem \ref{thm1.2} and Proposition \ref{pro1.1},
we can remove the smallness restriction (\ref{1.14}) on the initial
data and prove the following existence theorem of global strong solutions with large initial data, provided the
macroscopic molecular orientation of the liquid crystal materials
satisfies a natural geometric angle condition. This extends the Lei-Li-Zhang's result (cf.
\cite{LLZ2012}) to the case of nonhomogeneous incompressible liquid
crystal flows with initial vacuum.
\begin{thm}\label{thm1.3} Let $e_3=(0,0,1)\in \mathbb{S}^2$ and let $d_{03}$ be the third component of $d_0$. Besides
the condition \eqref{1.13} in Theorem \ref{thm1.1}, assume further that
\begin{equation}
 d_{03}\geq\varepsilon\quad {\rm and}\quad d_0-e_3\in L^2(\R^2)\label{1.20}
\end{equation}
holds for some uniform positive constant $\varepsilon>0$. Then for any $0<T<\infty$, there exists a unique
global strong solution $(\rho, u,P, d)$ of \eqref{1.1}--\eqref{1.4}, \eqref{1.10}
and \eqref{1.11} on $\R^2\times[0,T]$.
\end{thm}

The rest of the paper is organized as follows. In Sect. 2, we state some known inequalities and facts which will be used later.
The proof of Theorem \ref{thm1.1} will be done in Sect. 3, based on the local existence theorem and the global a priori estimates.
In Sect. 4, we outline the proof of Theorems \ref{thm1.2} and \ref{thm1.3}.

\section{Preliminaries}
In this section, we list some useful lemmas which will be frequently
used in the next sections. We first recall the well-known
Ladyzhenskaya and Sobolev inequalities (see, for example,
\cite{La1969,Ad1975}).
\begin{lem}\label{lem2.1} For $f\in H^1(\R^2)$, it holds for any $2\leq p<\infty$ that
\begin{eqnarray}
\|f\|_{L^4}^2&\leq& \sqrt 2\|f\|_{L^2}\|\nabla f\|_{L^2},\label{2.1}\\
\|f\|_{L^p}&\leq& C(p)\|f\|_{L^2}^{2/p}\|\nabla f\|_{L^2}^{1-2/p},\label{2.2}
\end{eqnarray}
where $C(p)$ is a positive constant depending on $p$. In addition, if $f\in W^{1,p}(\R^2)\cap H^2(\R^2)$ with $p>2$, then there exists a universal positive constant $C$ such that
\begin{equation}
\|f\|_{L^\infty}\leq C\|f\|_{W^{1,p}}\leq C\|f\|_{H^2}.\label{2.3}
\end{equation}
\end{lem}

We will also use the following Poincar${\rm\acute{e}}$ type
inequality, which shows that the velocity $u$ actually belongs to
$L^2$-space even that the vacuum states may appear.
\begin{lem}\label{lem2.2} Let $\tilde\rho>0$ be a given positive
constants. Assume that $\varrho-\tilde\rho\in L^2(\R^2)\cap
L^\infty(\R^2)$ with $\varrho(x)\geq0$, $\nabla v\in L^2(\R^2)$ and
$\sqrt\varrho v\in L^2(\R^2)$. Then,
\begin{equation}
\|v\|_{L^2}\leq C(\tilde\rho,\|\varrho-\tilde\rho\|_{{L^2}{\cap}{
L^\infty}})\left(\|\rho^{1/2}v\|_{L^2}+\|\nabla
v\|_{L^2}\right),\label{2.4}
\end{equation}
where $C(\tilde\rho,\|\varrho-\tilde\rho\|_{{L^2}{\cap}{
L^\infty}})$ is a positive constant depending only on $\tilde\rho$,
$\|\varrho-\tilde\rho\|_{L^2}$ and
$\|\varrho-\tilde\rho\|_{L^\infty}$.
\end{lem}
\pf Indeed, by virtue of H\"{o}lder and (\ref{2.2}), we have for any
$q\geq2$ that
\begin{eqnarray*}
\tilde\rho\int|v|^2dx&=&\int\varrho|v|^2dx-\int(\varrho-\tilde\rho)|v|^2dx\nonumber\\
&\leq&C\|\varrho^{1/2} v\|_{L^2}^2+C\left(\int|\varrho-\tilde\rho|^qdx\right)^{1/q}\left(\int|v|^{2q/(q-1)}dx\right)^{(q-1)/q}\nonumber\\
&\leq&C(\tilde\rho,\|\varrho-\tilde\rho\|_{{L^2}{\cap}{
L^\infty}})\left(\|\varrho^{1/2} v\|_{L^2}^2+\|v\|_{L^2}^{2(q-1)/q}\|\nabla v\|_{L^2}^{2/q}\right)\nonumber\\
&\leq&C(\tilde\rho,\|\varrho-\tilde\rho\|_{{L^2}{\cap}{
L^\infty}})\left(\|\varrho^{1/2} v\|_{L^2}^2+\|\nabla
v\|_{L^2}^2\right)+\frac{\tilde\rho}{2}\|v\|_{L^2}^2,
\end{eqnarray*}
which proves (\ref{2.4}) immediately.\hfill$\square$

\vskip 2mm

Next, to improve the regularity of the velocity, we need to use the
following estimates of the Stokes equations (see, for example,
\cite{Ga1994,LS2004}).
\begin{lem} \label{lem2.3}Consider the following stationary Stokes equations:
$$
-\Delta U+\nabla P=f,\quad
{\rm div} U=0\quad{in}\quad\R^2.
$$
Then for any $f\in W^{m,p}(\R^2)$ with $m\in \mathbb{Z}^+$ and $p>1$, there exists a positive constant $C$, depending only on $m$ and $p$, such that
\begin{equation}
\|\nabla^2U\|_{W^{m,p}}+\|\nabla P\|_{W^{m,p}}\leq
C\|f\|_{W^{m,p}}.\label{2.5}
\end{equation}
\end{lem}

To estimate the $L^2$-norm of the gradient of the velocity, we shall
apply a critical Sobolev inequality of logarithmic type which was
prove by Huang-Wang (cf. \cite{HW2013}) and is originally due to
Brezis-Wainger \cite{BW1980} (see also \cite{KOT2002,Oz1995}). This
is the key tool for the proofs of Theorems \ref{thm1.1}--\ref{1.3}.
\begin{lem}\label{lem2.4} For $q>2$ and $0\leq s<t<\infty$, assume that $f\in L^2(s,t;H^1(\R^2))\cap L^2(s,t;W^{1,q}(\R^2))$. Then there exists a positive constant $C(q)$, independent of $s,t$, such that
\begin{equation}
\|f\|_{L^2(s,t;L^\infty(\R^2))}\leq
C\left(1+\|f\|_{L^2(s,t;H^1(\R^2))}\left(\ln
^+\|f\|_{L^2(s,t;W^{1,q}(\R^2))}\right)^{1/2}\right).\label{2.6}
\end{equation}
\end{lem}

In the case that the lower bound of the density is nonnegative, the
local existence of strong solutions to (\ref{1.1})--(\ref{1.4}),
(\ref{1.10}) and (\ref{1.11}) was proved in  \cite{WD2011}. Indeed,
in \cite{WD2011} the authors only considered the case of smooth
bounded domains, however, as pointed out in \cite{CK2003}, the
similar procedure also works for the whole space by means of the
standard domain expansion technique. For simplicity, we quote the
following local existence theorem of strong solutions without
proofs.
\begin{lem}\label{lem2.5}Assume that the conditions of Theorem \ref{thm1.1} hold. Then there exists a positive time $0<T_0<\infty$ such that the Cauchy problem \eqref{1.1}--\eqref{1.4}, \eqref{1.10}
and \eqref{1.11} admits a unique strong solution on $\R^2\times(0,T_0)$.

\end{lem}

\section{Proof of Theorem \ref{thm1.1}}
Assume that the conditions of Theorem \ref{thm1.1} hold. Let $0<T^*<\infty$ be the first blowup time of a
strong solution $(\rho,u,P,d)$ to the Cauchy problem
(\ref{1.1})--(\ref{1.4}), (\ref{1.10}) and (\ref{1.11}). In order to prove Theorem \ref{thm1.1}, it suffices to prove there actually exists a generic positive
constant $0<M<\infty$, depending only on the initial data
$(\rho_0,u_0,d_0)$ and $T^*$, such that
\begin{eqnarray}
\mathcal {E}(T)&\triangleq&\sup_{0\leq t\leq
T}\left(\|\rho-\tilde\rho\|_{H^2}+\|u\|_{H^2}+\|\nabla
d\|_{H^2}+\|\rho^{1/2}
u_t\|_{L^2}^2+\|d_t\|_{H^1}^2\right)\nonumber\\
&&\quad+\int_0^T\left(\|u\|_{H^3}^2+\|\nabla
d\|_{H^3}^2+\|u_t\|_{H^1}^2+\|d_t\|_{H^2}^2\right)dt\leq
M\label{3.1}
\end{eqnarray}
holds for any $0<T<T^*$. So, by the local existence theorem (see Lemma \ref{lem2.5}) it can
be easily shown that the strong solution can be extended beyond
$T^*$, which gives a contradiction of $T^*$. Hence, the strong solution
exists globally on $\R^2\times[0,T]$ for any $0<T<\infty$. The proof of Theorem \ref{thm1.1} is therefore complete.

The proof of (\ref{3.1}) is based on a series of lemmas.
Throughout the remainder of the paper, for simplicity we denote by $C$ a generic
constant which depends only on the initial data and $T^*$,
and may change from line to line.

First, it is easy to see from the method of characteristics
and (\ref{1.1}) that for every $0<T<T^*$,
\begin{equation}
0\leq\rho(x,t)\leq \|\rho_0\|_{L^\infty}\quad{\rm for\;\;all}\quad
(x,t)\in\R^2\times[0,T].\label{3.2}
\end{equation}
Moreover, multiplying (\ref{1.1}) by
$q|\rho-\tilde\rho|^{q-2}(\rho-\tilde\rho)$ with $q\geq2$,
integrating it by parts over $(0,t)$, and using the divergence-free
condition (\ref{1.3}), we find that
\begin{equation}
\|(\rho-\tilde\rho)(t)\|_{L^q}=\|\rho_0-\tilde\rho\|_{L^q}\quad{\rm
for}\quad\forall\; t\in[0,T].\label{3.3}
\end{equation}

In view of (\ref{1.1})--(\ref{1.4}), we have the following standard
energy estimates.
\begin{lem}\label{lem3.1} For every $0<T <T^*$, one has
\begin{eqnarray}
&&\sup_{0\leq t\leq
T}\int \left(|\rho^{1/2}u|^2+ |\nabla d|^2\right)dx+2\int_0^T\int \left(|\nabla u|^2+|\Delta d+|\nabla d|^2d|^2\right) d x d t \nonumber\\
&&\qquad\leq \int \left(|\rho^{1/2}u|^2+ |\nabla
d|^2\right)(x,0)dx\triangleq E_0.\label{3.4}
\end{eqnarray}
\end{lem}
\pf Multiplying $\eqref{1.2}$ by $u$ in $L^2$ and integrating by
parts, by (\ref{1.3}) we know that
\begin{eqnarray}
\frac{1}{2}\frac{d }{d t}\int  \rho|u|^2 d x+ \int_\Omega |\nabla
u|^2  d x=-\int  \left(u\cdot\nabla d \cdot\Delta d\right)  d
x.\label{3.5}
\end{eqnarray}
Due to the fact that $|d|=1$, multiplying \eqref{1.4}  by $(\Delta
d+|\nabla d|^2d)$ in $L^2$, we obtain after integrating the
resulting equations by parts over $\R^2$ that
\begin{eqnarray}
&&\frac{1}{2}\frac{d }{ d t}\int  |\nabla d|^2  d x+\int_\Omega |\Delta d+|\nabla d|^2d|^2  d x\nonumber\\
&&\quad=\int \left( u\cdot\nabla d\cdot\Delta d\right)  d x
+\int \left(|\nabla d|^2d\cdot d_t+|\nabla d|^2u\cdot\nabla d\cdot d \right) d x\nonumber\\
&&\quad=\int \left( u\cdot\nabla d\cdot\Delta d\right)  d x
+\frac{1}{2}\int \left(|\nabla d|^2\partial_t|d|^2 +|\nabla d|^2u\cdot\nabla |d|^2 \right) d x\nonumber\\
&&\quad=\int \left( u\cdot\nabla d\cdot\Delta d\right)  d
x,\label{3.6}
\end{eqnarray}
which, combined with (\ref{3.5}), immediately leads to
\eqref{3.4}.\hfill$\square$

\vskip 2mm

To be continued, we need the following key estimates on $\|\nabla^2
d\|_{L^2(0,T;L^2)}$.
\begin{lem}\label{lem3.2}Assume that the initial data satisfies
\begin{eqnarray}
 \exp{\left(2\left(\|\rho_0^{1/2}u_0\|^2_{L^2}+\|\nabla d_0\|^2_{L^2}\right)\right)}\|\nabla d_0\|^2_{L^2}\leq \frac{1}{16},\label{3.7}
\end{eqnarray}
then it holds for every $T\in(0,T^*)$ that
\begin{equation} \label{3.8}
\sup_{0\leq t\leq T}\|\nabla d\|^2_{L^2}+ \int_0^T\|\nabla^2
d\|^2_{L^2} dt\leq \frac{1}{16}.
\end{equation}
\end{lem}
\pf After integrating by parts, we easily deduce from the identity
$|d|=1$ that
\begin{eqnarray}
\int|\Delta d+|\nabla d|^2d|^2dx&=&\int\left(|\Delta d|^2+|\nabla
d|^4\right) dx-2\int|\nabla d|^2 (d\cdot\Delta
d)dx\nonumber\\
&=&\int\left(|\Delta d|^2-|\nabla d|^4\right) dx.\label{3.9}
\end{eqnarray}
On the other hand, integration by parts, together with the
divergence-free condition (\ref{1.3}), gives
\begin{eqnarray}
\int \left( u\cdot\nabla d\cdot\Delta d\right)  d x&=&-\int \left(\partial_j u^i\partial_i d^k\partial_j d^k+u^i\partial_{ij}^2 d^k\partial_j d^k\right)  d x\nonumber\\
&=&-\int\left(\partial_j u^i\partial_i d^k\partial_j d^k\right)  d
x\leq\|\nabla u\|_{L^2}\|\nabla d\|_{L^4}^2,\label{3.10}
\end{eqnarray}
where and in what follows the repeated indices denotes the summation
over the indices.

Putting (\ref{3.9}), (\ref{3.10}) into (\ref{3.6}) and recalling the
fact that
$$
\|\Delta d\|_{L^2}^2=\|\nabla^2 d\|_{L^2}^2,
$$
we obtain
$$
\frac{1}{2}\frac{d}{dt}\|\nabla d\|_{L^2}^2+\|\nabla^2
d\|_{L^2}^2\leq \|\nabla u\|_{L^2}\|\nabla d\|_{L^4}^2+\|\nabla
d\|_{L^4}^4,
$$
which, combined with  (\ref{2.1}) and the Cauchy-Schwarz inequality,
yields
\begin{eqnarray}
\frac{1}{2}\frac{d}{dt}\|\nabla d\|_{L^2}^2+\|\nabla^2 d\|_{L^2}^2&\leq& \sqrt 2\|\nabla u\|_{L^2}\|\nabla d\|_{L^2}\|\nabla^2d\|_{L^2}+2\|\nabla d\|_{L^2}^2\|\nabla^2 d\|_{L^2}^2\nonumber\\
&\leq&\left(2\|\nabla
d\|_{L^2}^2+\frac{1}{4}\right)\|\nabla^2d\|_{L^2}^2+2\|\nabla
u\|_{L^2}^2\|\nabla d\|_{L^2}^2.\label{3.11}
\end{eqnarray}

It follows from (\ref{3.7}) that
$$
\|\nabla d_0\|_{L^2}^2\leq e^{2E_0}\|\nabla d_0\|_{L^2}^2\leq
\frac{1}{16},
$$
and thus, by the local existence theorem and the continuity argument
we see that there exists a $T_1>0$ such that for any $t\in[0,T_1]$,
\begin{equation}
\|\nabla d\|_{L^2}^2\leq \frac{1}{8}.\label{3.12}
\end{equation}

Set
$$
\tilde T\triangleq\sup_{}\{ T\; |\; (\ref{3.12})\; \; {\rm holds}\}.
$$
Then it follows from (\ref{3.11})--(\ref{3.12}) that for any
$t\in[0,\tilde T)$,
$$
\frac{d}{dt}\|\nabla d\|_{L^2}^2+\|\nabla^2 d\|_{L^2}^2\leq
4\|\nabla u\|_{L^2}^2\|\nabla d\|_{L^2}^2,
$$
which, together with Gronwall's inequality and (\ref{3.4}), leads to
\begin{eqnarray}
\|\nabla d\|_{L^2}^2+ \int_0^t\|\nabla^2 d\|_{L^2}^2 d\tau&\leq&  \exp\left(4\int_0^t\|\nabla u\|_{L^2}^2d\tau \right)\|\nabla d_0\|_{L^2}^2\nonumber\\
&\leq& e^{2E_0}\|\nabla d_0\|_{L^2}^2\leq \frac{1}{16}.\label{3.13}
\end{eqnarray}
Combining (\ref{3.4}), (\ref{3.13}) with the continuity argument
immediately implies that (\ref{3.12}) holds for all $0<\tilde
T<T^*$, and thus, the proof of (\ref{3.8}) is finished.\hfill
$\square$

\vskip 2mm

By Lemmas \ref{lem2.4} and \ref{lem3.2}, we can now derive the
estimates of $\|\nabla u\|_{L^2}$ and $\|\nabla^2 d\|_{L^2}$ which
is the most important step among the proofs.
\begin{lem}\label{lem3.3} For every $0<T<T^*$, one has
\begin{eqnarray}
&&\sup_{0< t\leq T}\left(\|u\|_{H^1}^2+\|\nabla
d\|_{H^1}^2+\|d_t\|_{L^2}^2\right)\nonumber\\
&&\qquad+\int_0^T\left(\|\nabla^2u\|_{L^2}^2+\|\rho^{1/2}\dot
u\|_{L^2}^2+\|d_t\|_{H^1}^2+\|\nabla  d\|_{H^2}^2\right) d t \leq
C,\label{3.14}
\end{eqnarray}
which particularly gives
\begin{equation}
\int_0^T\|\rho^{1/2}u_t\|_{L^2}^2\leq C.\label{3.15}
\end{equation}
\end{lem}
\pf Let $\dot f\triangleq f_t+u\cdot\nabla f$ denote the material derivative. Also set
$$
M(d)\triangleq \nabla d\otimes \nabla d-\frac{1}{2}|\nabla d|^2\mathbb{I}_2,\quad (\nabla d\otimes \nabla d)_{ij}\triangleq\frac{\partial d}{\partial x_i}\cdot\frac{\partial d}{\partial x_j},\quad 1\leq i,j\leq 2,$$
then it is easily seen that
$$\nabla d\cdot\Delta d={\rm div}(M(d)).$$

To prove (\ref{3.14}), multiplying $(\ref{1.2})$ by $u_t$ and integrating it by parts over
$\R^2$, we deduce
\begin{eqnarray*}
&&\frac{1}{2}\frac{ d}{ d t}\|\nabla
u\|_{L^2}^2+\|\rho^{1/2}\dot u\|_{L^2}^2
\nonumber\\
&&\quad=-\int {\rm div}M(d)\cdot u_t d x+\int\rho u\cdot\nabla u\cdot
\dot u d x\nonumber\\
&&\quad=\frac{d}{dt}\int  M(d):\nabla u d x-\int M(d)_t:\nabla u d x+\int\rho u\cdot\nabla u\cdot
\dot u d x\nonumber\\
&&\quad\leq\frac{d}{dt}\int  M(d):\nabla u d x+\frac{1}{2}\|\rho^{1/2}\dot u\|^2_{L^2}
+\frac{1}{4}\|\nabla d_t\|^2_{L^2}+C(\|u\|^2_{L^\infty}+\|\nabla d\|^2_{L^\infty})\|\nabla u\|^2_{L^2},
\end{eqnarray*}
where we have also used (\ref{3.2}) and Cauchy-Schwarz inequality.
As a result,
\begin{eqnarray}
&& \frac{ d}{ d t}\|\nabla
u\|_{L^2}^2+ \|\rho^{1/2}\dot u\|_{L^2}^2
\nonumber\\
&&\quad\leq 2\frac{d}{dt}\int  M(d):\nabla u d x+\frac{1}{2}\|\nabla
d_t\|^2_{L^2}+C(\|u\|^2_{L^\infty} +\|\nabla
d\|^2_{L^\infty})\|\nabla u\|^2_{L^2}.\label{3.16}
\end{eqnarray}

Next, one easily obtains from (\ref{1.4}) that
\begin{eqnarray}
\frac{ d}{ d t}\|\nabla
d\|_{L^2}^2+\left(\|d_t\|_{L^2}^2+\|\nabla^2d\|_{L^2}^2\right)&\leq&
C\int\left(|u|^2 |\nabla d|^2+|\nabla
d|^4\right)dx\nonumber\\
&\leq& C\left(\|u\|_{L^\infty}^2+\|\nabla
d\|_{L^\infty}^2\right)\|\nabla d\|_{L^2}^2.\label{3.17}
\end{eqnarray}

To deal with the term $\|\nabla d_t\|^2_{L^2}$ on the right-hand
side of \eqref{3.16}, we first apply $\nabla$ to  both sides of
\eqref{1.4} to get that
\begin{equation}
\nabla d_t-\nabla\Delta d=-\nabla(u\cdot\nabla d)+\nabla(|\nabla
d|^2d),\label{3.18}
\end{equation}
from which it follows that
\begin{eqnarray}
&&\frac{ d}{ d t}\|\nabla^2
d\|_{L^2}^2+\left(\|\nabla  d_t\|_{L^2}+\|\nabla \Delta d\|_{L^2}^2\right)\nonumber\\
&&\quad\leq \int \left(|\nabla(|\nabla d|^2 d)|^2+|\nabla (u\cdot \nabla d)|^2\right) d x\nonumber\\
&&\quad\leq C\int\left(|\nabla d|^6+|\nabla d|^2|\nabla^2 d|^2+|\nabla u|^2|\nabla d|^2 +|u|^2|\nabla^2 d|^2\right) dx\nonumber\\
&&\quad\leq C\left(\|u\|^2_{L^\infty}+\|\nabla
d\|^2_{L^\infty}\right)\left(\|\nabla u\|_{L^2}^2+\|\nabla^2
d\|_{L^2}^2\right),\label{3.19}
\end{eqnarray}
where we have used (\ref{2.1}) and (\ref{3.4}) to get that
$$
\|\nabla d\|_{L^6}^6\leq C\|\nabla d\|_{L^\infty}^2\|\nabla d\|_{L^2}^2\|\nabla^2 d\|_{L^2}^2
\leq C\|\nabla d\|_{L^\infty}^2\|\nabla^2 d\|_{L^2}^2.
$$

Using (\ref{2.1}) and (\ref{3.4}) again, we have
$$
\int  M(d):\nabla u d x\leq \frac{1}{4}\|\nabla u\|^2_{L^2}+C\|\nabla d\|^2_{L^2}\|\nabla^2 d\|^2_{L^2}
\leq \frac{1}{4}\|\nabla u\|^2_{L^2}+C_1\|\nabla^2 d\|^2_{L^2}.
$$
Taking this into account,  multiplying \eqref{3.19} by $2C_1 +1$,
and adding the resulting inequality, (\ref{3.16}) and \eqref{3.17}
together, we obtain after integrating the resulting inequality over
$(s,t)$ with $0\leq s<t<T$ that
\begin{eqnarray*}
&&\left(\|\nabla u\|_{L^2}^2 +\|\nabla
d\|_{H^1}^2\right)(t)+\int_s^t\left(\|\rho^{1/2}\dot
u\|_{L^2}^2+\|d_t\|_{H^1}^2+\|\nabla^2 d\|_{H^1}^2\right)d
\tau\nonumber
\\&&\quad\leq C\left(\|\nabla
u\|_{L^2}^2+\|\nabla
d\|_{H^1}^2\right)(s)+C\int_s^t\left(\|u\|^2_{L^\infty}+\|\nabla
d\|^2_{L^\infty}\right)\left(\|\nabla u\|^2_{L^2}+\|\nabla
d\|_{H^1}^2\right) d \tau ,
\end{eqnarray*}
and consequently,
\begin{eqnarray}
&&\left(\|\nabla u\|_{L^2}^2 +\|\nabla
d\|_{H^1}^2\right)(t)+\int_s^t\left(\|\rho^{1/2}\dot
u\|_{L^2}^2+\|d_t\|_{H^1}^2+\|\nabla^2 d\|_{H^1}^2\right)d
\tau\nonumber
\\&&\quad\leq C\left(\|\nabla
u\|_{L^2}^2+\|\nabla
d\|_{H^1}^2\right)(s)\exp\left(C\int_s^t\left(\|u\|^2_{L^\infty}+\|\nabla
d\|^2_{L^\infty}\right) d \tau \right).\label{3.20}
\end{eqnarray}

Clearly, it remains to estimate $\|(u,\nabla d)\|_{L^\infty}$. To
this end, let
$$
\Phi(t) \triangleq e+\sup_{0\leq\tau\leq t}\left(\|\nabla
u\|_{L^2}^2 +\|\nabla
d\|_{H^1}^2\right)(\tau)+\int_0^t\left(\|\rho^{1/2}\dot
u\|_{L^2}^2+\|d_t\|_{H^1}^2+\|\nabla^2 d\|_{H^1}^2\right)d \tau.
$$

First, in view of (\ref{2.4}) and (\ref{3.2})--(\ref{3.4}), we have
\begin{equation}
\|u\|_{L^2}^2\leq C\left(\|\rho^{1/2}u\|_{L^2}^2+\|\nabla
u\|_{L^2}^2\right)\leq  C\left(1+\|\nabla
u\|_{L^2}^2\right).\label{3.21}
\end{equation}

Next, using Lemma \ref{lem2.3}, (\ref{2.1}), (\ref{3.2}) and
(\ref{3.4}), we deduce from H\"{older} and Cauchy-Schwarz
inequalities that
\begin{eqnarray}\label{3.22}
\|\nabla^2u\|_{L^{2}}&\leq& C(\|\rho \dot u\|_{L^{2}}+\|\nabla d\cdot\Delta d\|_{L^{2}})\nonumber\\
&\leq& C\left(\|\rho^{1/2} \dot u\|_{L^2}+\|\nabla^2
d\|_{L^2}^2+\|\nabla^3 d\|_{L^2}\right),
\end{eqnarray}
which, combined with (\ref{3.8}), yields
\begin{eqnarray}\label{3.23}
\int_s^t\|\nabla^2u\|_{L^2}^2d\tau&\leq& C\int_s^t\left(\|\rho^{1/2}
\dot u\|_{L^2}^2+\|\nabla^2 d\|_{L^2}^4+\|\nabla^3
d\|_{L^2}^2\right)\dd \tau\nonumber\\
&\leq&C\sup_{s\leq \tau\leq t} \|\nabla^2
d\|_{L^2}^2+C\int_s^t\left(\|\rho^{1/2} \dot u\|_{L^2}^2+\|\nabla^3
d\|_{L^2}^2\right)\dd \tau.
\end{eqnarray}

Thus, recalling the definition of $\Phi(T)$ and using (\ref{3.4}),
(\ref{3.21}) and (\ref{3.23}), we infer from Lemma \ref{lem2.4} that
for any $0\leq s<t\leq T<T^*$,
\begin{eqnarray}
\|u\|_{L^2(s,t;L^\infty)}^2&\leq& C\left(1+\|u\|_{L^2(s,t;H^1)}^2\ln\left(e+\|u\|_{L^2(s,t;W^{1,4})}\right)\right)\nonumber\\
&\leq& C \left(1+\|\nabla
u\|_{L^2(s,t;L^2)}^2\ln \left(e+\|u\|_{L^2(s,t;H^2)}\right)\right)\nonumber\\
&\leq& C \left(1+\|\nabla u\|_{L^2(s,t;L^2)}^2\ln
\left(e+\|\nabla^2 u\|_{L^2(s,t;L^2)}\right)\right)\nonumber\\
&\leq& C \left(1+\|\nabla u\|_{L^2(s,t;L^2)}^2\ln \left(C
\Phi(t)\right)\right).\label{3.24}
\end{eqnarray}
In a similar manner, by (\ref{2.6}) and (\ref{3.4}) one has
\begin{eqnarray}
\|\nabla d\|_{L^2(s,t;L^\infty)}^2&\leq& C\left(1+\|\nabla d\|_{L^2(s,t;H^1)}^2\ln\left(e+\|\nabla d\|_{L^2(s,t;W^{1,4})}\right)\right)\nonumber\\
&\leq& C\left(1+\|\nabla^2
d\|_{L^2(s,t;L^2)}^2\ln \left(e+\|\nabla^2 d\|_{L^2(s,t;H^1)}\right)\right)\nonumber\\
&\leq& C \left(1+\|\nabla^2 d\|_{L^2(s,t;L^2)}^2\ln \left(C
\Phi(t)\right)\right).\label{3.25}
\end{eqnarray}

For any $0\leq s<t\leq T<T^*$, putting (\ref{3.24}) and (\ref{3.25})
into (\ref{3.20}) gives
\begin{eqnarray}
\Phi(t)&\leq& C \Phi(s)\exp\left\{C_2 \left(\|\nabla
u\|_{L^2(s,t;L^2)}^2+\|\nabla^2 d\|_{L^2(s,t;L^2)}^2\right)\ln
\left(C_1 \Phi(t)\right)\right\}\nonumber\\
&\leq&C \Phi(s)\left[ C_1 \Phi(t) \right]^{C_2 \left(\|\nabla
u\|_{L^2(s,t;L^2)}^2+\|\nabla^2
d\|_{L^2(s,t;L^2)}^2\right)}.\label{3.26}
\end{eqnarray}

It follows from (\ref{3.4})  and (\ref{3.8}) that there exists a
positive constant $\delta>0$ such that
$$
C_2 \left(\|\nabla u\|_{L^2(T-\delta,T;L^2)}^2+\|\nabla^2
d\|_{L^2(T-\delta,T;L^2)}^2\right)\leq\frac{1}{2},
$$
which, inserted into (\ref{3.26}), leads to
$$
\Phi(T)\leq C \Phi(T-\delta)\left[C_1 \Phi(T)\right]^{1/2}\leq
\frac{1}{2}\Phi(T)+C \Phi^2(T-\delta),
$$
so that
\begin{equation}
\Phi(T)\leq C(T)\Phi^2(T-\delta).\label{3.27}
\end{equation}

As a result of (\ref{3.27}), we see that $\Phi(T)$ is bounded for
any $0<T<T^*$ since  the local existence theorem indicates
$\Phi(T-\delta)<\infty $ for any $0<T<T^*$. This, together with
(\ref{3.21}) and (\ref{3.22}), finishes the proof of (\ref{3.14}).

Furthermore, recalling the definition of material derivative (i.e.
`` $\dot{}$ ''), one gets from  (\ref{2.1}), (\ref{3.2}) and
(\ref{3.14}) that
\begin{eqnarray*}
\int_0^T\|\rho^{1/2}u_t\|_{L^2}^2dt&\leq& \int_0^T\|\rho^{1/2}\dot u\|_{L^2}^2dt+\int_0^T\int\rho|u|^2|\nabla u|^2dx dt\\
&\leq& C+C\int_0^T\|u\|_{L^4}^2\|\nabla u\|_{L^4}^2dt\\
&\leq&C+C\int_0^T\|u\|_{L^2}\|\nabla u\|_{L^2}^2\|\nabla^2u\|_{L^2}dt\\
&\leq&C+C\int_0^T \|\nabla^2u\|_{L^2}^2dt\leq C,
\end{eqnarray*}
which immediately proves (\ref{3.15}). The proof of Lemma
\ref{lem3.3} is therefore complete.\hfill$\square$

\vskip 2mm

Next, we proceed to estimate $\|\rho^{1/2}u_t\|_{L^2}$ and $\|\nabla
d_t\|_{L^2}$.
\begin{lem}\label{lem3.4} For every $0<T<T^*$, one has
\begin{eqnarray}
\sup_{0< t\leq
T}\left(\|\rho^{1/2}u_t\|_{L^2}^2+\|
d_t\|_{H^1}^2\right)+\int_0^T\left(\|\nabla
u_t\|_{L^2}^2+\|d_{tt}\|_{L^2}^2+\|\nabla^2d_t\|_{L^2}^2\right) dt \leq C,\label{3.28}
\end{eqnarray}
and moreover,
\begin{eqnarray}
\sup_{0< t\leq
T}\left(\|u\|_{H^{2}}^2+\|\nabla d\|_{H^{2}}^2\right)+\int_0^T\left(\|\nabla u\|^2_{W^{1,4}}+\|\nabla^2d\|_{H^2}^2\right)dt\leq C.\label{3.29}
\end{eqnarray}
\end{lem}
\pf Differentiating \eqref{1.2} with respect to $t$ gives
$$
\rho u_{tt}+\rho u\cdot\nabla u_t-\Delta u_t=-\rho_{t}
(u_{t}+u\cdot\nabla u)-\rho u_t\cdot\nabla u-\nabla P_t-{\rm div}M_t,
$$
which, multiplied by $u_t$ in $L^2$ and integrated by parts over $\R^2$, results in
\begin{eqnarray}
&&\frac{1}{2}\frac{ d}{ d t}\|\rho^{1/2}u_t\|_{L^2}^2+\|\nabla u_t\|_{L^2}^2\nonumber\\
&&\quad=-\int\rho_t(u_t+u\cdot\nabla u)\cdot u_t d x-\int\rho u_t\cdot\nabla u\cdot u_t dx+\int M_t:\nabla u_t d x\nonumber\\
&&\quad\triangleq I_1+I_2+I_3\label{3.30}
\end{eqnarray}
where $A:B=\sum_{i,j=1}^2a_{ij}b_{ij}$ for $A=(a_{ij})_{2\times 2}$
and $B=(b_{ij})_{2\times2}$.

We are now in a position of estimating the right-hand side of (\ref{3.30}) term by term.
First, using (\ref{1.1}) and integrating by parts, by Lemma \ref{lem2.1}, (\ref{3.2}) and (\ref{3.14}) we deduce
\begin{eqnarray*}
I_1&=& \int\left(\rho u\cdot\nabla |u_t|^2 +\rho u \cdot\nabla(u\cdot\nabla u\cdot u_t)\right) d x \nonumber\\
&\leq& C\int\left(\rho|u||u_t||\nabla u_t|+\rho|u||\nabla u|^2|u_t|+\rho|u|^2|\nabla^2 u||u_t|+\rho|u|^2|\nabla u||\nabla u_t|\right)dx\nonumber\\
&\leq&C\left(\|u\|_{L^\infty}\|\rho^{1/2}u_t\|_{L^2}\|\nabla u_t\|_{L^2}+\|u\|_{L^\infty}\|\nabla u\|_{L^4}^2\|\rho^{1/2} u_t\|_{L^2}\right)\nonumber\\
&&+C\left(\|u\|_{L^\infty}^2\|\rho^{1/2} u_t\|_{L^2}\|\nabla^2u\|_{L^2}+\|u\|_{L^\infty}^2\|\nabla u\|_{L^2}\|\nabla u_t\|_{L^2}\right)\nonumber\\
&\leq&\frac{1}{4}\|\nabla u_t\|_{L^2}^2+C\left(1+\|\rho^{1/2}u_t\|_{L^2}^4+\|\nabla^2 u\|_{L^2}^2\right),
\end{eqnarray*}
where we have used Cauchy-Schwarz inequality and the following estimate due to Lemma \ref{lem2.1} and (\ref{3.14}):
\begin{equation}
\|u\|_{L^\infty}\leq C\|u\|_{W^{1,4}}\leq C\left(\|u\|_{H^1}+\|\nabla u\|_{L^2}^{1/2}\|\nabla^2u\|_{L^2}^{1/2}\right)\leq C\left(1+\|\nabla^2u\|_{L^2}^{1/2}\right).\label{3.31}
\end{equation}

Due to (\ref{2.4}), (\ref{3.2}) and (\ref{3.3}), we have
\begin{equation}\label{3.32}
\|u_t\|_{L^2}^2\leq C\left(\|\rho^{1/2} u_t\|_{L^2}^2+\|\nabla u_t\|_{L^2}^2\right),
\end{equation}
and thus, by (\ref{2.1}), (\ref{3.2}) and (\ref{3.14}) the second term $I_2$ can be bounded as follows:
\begin{eqnarray*}
|I_2|&\leq& C\int\rho |u_t|^2|\nabla u|\dd x\nonumber\\
&\leq& C\|\nabla u\|_{L^4}\|u_t\|_{L^4}\|\rho^{1/2}u_t\|_{L^2}\nonumber\\
&\leq& C\|\nabla u\|_{L^2}^{1/2}\|\nabla^2u\|_{L^2}^{1/2}\|u_t\|_{L^2}^{1/2}\|\nabla u_t\|_{L^2}^{1/2}\|\rho^{1/2}u_t\|_{L^2}\nonumber\\
&\leq& C\left(\|\nabla u_t\|_{L^2}+\|\rho^{1/2}u_t\|_{L^2}\right)\|\nabla^2 u\|_{L^2}^{1/2}\|\rho^{1/2}u_t\|_{L^2}\nonumber\\
&\leq& \frac{1}{4}\|\nabla u_t\|^2_{L^2}+C\left(1+\|\nabla^2 u\|^2_{L^2}+\|\rho^{1/2}u_t\|^4_{L^2}\right).
\end{eqnarray*}

Finally, it is easily seen from (\ref{2.3}) and (\ref{3.14}) that
\begin{eqnarray*}
|I_3|&\leq& C\|\nabla d\|_{L^\infty}\|\nabla d_t\|_{L^2}\|\nabla u_t\|_{L^2}\nonumber\\
&&\leq \frac
{1}{4}\|\nabla u_t\|^2_{L^2}+C\left(1+\|\nabla^3 d\|^2_{L^2}\right)\|\nabla d_t\|^2_{L^2}.
\end{eqnarray*}

Substituting the estimates of $I_1,I_2$ and $I_3$ into \eqref{3.30}, one obtains
\begin{eqnarray}
&&\frac{d}{d t}\|\rho^{1/2}u_t\|_{L^2}^2+\|\nabla u_t\|_{L^2}^2\nonumber\\
&&\quad \leq C\left(1+\|\nabla^2u\|_{L^2}^2+\|\rho^{1/2} u_t\|^4_{L^2}\right)+C\left(1+\|\nabla d\|^2_{H^2}\right)\|\nabla d_t\|^2_{L^2}.\label{3.33}
\end{eqnarray}

To estimate $\|\nabla d_t\|_{L^2}$, we
differentiate (\ref{1.4}) with respect to $t$ to get
$$
d_{tt}-\Delta d_t=\left(|\nabla d|^2d-u\cdot\nabla d\right)_t,
$$
and hence, using Lemma \ref{lem2.1} and (\ref{3.14}), we deduce after direct calculations that
\begin{eqnarray}
&& \frac{ d}{ d t}\|\nabla d_t\|_{L^2}^2+\left(\|d_{tt}\|_{L^2}^2+\|\Delta d_t\|_{L^2}^2\right)\nonumber\\
&&\quad\leq C\int \left(|\nabla d|^2|\nabla d_t|^2+|\nabla d|^4|d_t|^2+|u_t|^2|\nabla d|^2+|u|^2|\nabla d_t|^2\right) d x\nonumber\\
&&\quad\leq C\left(\|\nabla
d\|_{L^\infty}^2+\|u\|_{L^\infty}^2\right)\|\nabla
d_t\|_{L^2}^2+C\|\nabla d\|_{L^\infty}^4\|d_t\|_{L^2}^2+C\|\nabla
d\|_{L^4}^2\|u_t\|_{L^4}^2\nonumber\\
&&\quad\leq
C\left( \|\nabla d\|_{H^2}^2+\|u\|_{H^2}^2\right)\|\nabla
d_t\|_{L^2}^2+C_1\left(\|\nabla d\|_{H^2}^2+\|\rho^{1/2}u_t\|_{L^2}^2+\|\nabla
u_t\|_{L^2}^2\right),\label{3.34}
\end{eqnarray}
where we have also used (\ref{3.32}) and the following estimate due to (\ref{2.1}), (\ref{2.3}) and (\ref{3.14}):
$$
\|\nabla d\|_{L^\infty}^4\leq C\left(\|\nabla
d\|_{L^4}^4+\|\nabla^2d\|_{L^4}^4\right)\leq C\left(1+\|\nabla^2
d\|_{L^2}^2\|\nabla^3 d\|_{L^2}^2\right)\leq
C\left(1+\|\nabla^3d\|_{L^2}^2\right).
$$

Now, multiplying \eqref{3.33} by $2C_1+1$ and adding it to \eqref{3.34}, we see that
\begin{eqnarray*}
&&\frac{ d}{d t}\left(\|\rho^{1/2}u_t\|_{L^2}^2+\|\nabla d_t\|_{L^2}^2\right)+ \left(\|\nabla u_t\|_{L^2}^2+\|\Delta d_t\|_{L^2}^2+\|d_{tt}\|_{L^2}^2\right)\nonumber\\
&&\quad\leq
C\left(1+\|\nabla d\|_{H^2}^2+\|u\|_{H^2}^2\right)\left(1+\|\nabla
d_t\|_{L^2}^2\right)+C\|\rho^{1/2}u_t\|_{L^2}^4 ,
\end{eqnarray*}
which, combined with \eqref{3.14},
\eqref{3.15} and Gronwall's inequality, leads to \eqref{3.28}, since the compatibility condition stated in (\ref{1.13})$_2$ implies that $(\rho^{1/2}u_t)(x,0)\in L^2(\R^2)$ is well defined.

Using (\ref{2.1}), (\ref{3.2}), (\ref{3.14}) and (\ref{3.28}),  we have by (\ref{2.5}) that
\begin{eqnarray}
\|\nabla^2 u\|_{L^2}&\leq& C\left(\|\rho u_t\|_{L^2}+\|\rho u\cdot\nabla u\|_{L^2}+\|\nabla d\cdot\Delta d\|_{L^2}\right)\nonumber\\
&\leq& C\left(\|\rho^{1/2}u_t\|_{L^2}+\|u\|_{L^4}\|\nabla u\|_{L^4}+\|\nabla d\|_{L^4}\|\nabla^2d\|_{L^4}\right)\nonumber\\
&\leq&C\left(1+\|\nabla^2 u\|_{L^2}^{1/2}+\|\nabla^3 d\|_{L^2}^{1/2}\right).\label{3.35}
\end{eqnarray}
Similarly, one also infers from (\ref{3.18}) that
\begin{eqnarray*}
\|\nabla^3 d\|_{L^2}&\leq&C\left(\|\nabla d_t\|_{L^2}+\|\nabla (u\cdot \nabla d)\|_{L^2}+\|\nabla (|\nabla d|^2d)\|_{L^2}\right)\nonumber\\
&\leq& C\left(1+\|u\|_{L^4}\|\nabla^2 d\|_{L^4}+\|\nabla
u\|_{L^4}\|\nabla d\|_{L^4}+\|\nabla d\|_{L^6}^3+\|\nabla d\|_{L^4}\|\nabla^2d\|_{L^4}\right)\nonumber\\
&\leq& C\left(1+\|\nabla^2u\|_{L^2}^{1/2}+\|\nabla^3d\|_{L^2}^{1/2}\right),
\end{eqnarray*}
from which, (\ref{3.35}) and Young's inequality, we arrive at
\begin{equation}
\sup_{0\leq t\leq T}\left(\|u\|_{H^2}+\|\nabla d\|_{H^2}\right)\leq C.\label{3.36}
\end{equation}

Using (\ref{2.5}) and Lemma \ref{lem2.1} again, we deduce from
(\ref{3.2}), (\ref{3.28}) and (\ref{3.36}) that
\begin{eqnarray}
\int_0^T\|\nabla^2 u\|_{L^4}^2&\leq&C\int_0^T\left(\|\rho u_t\|_{L^4}^2+\|\rho u\cdot\nabla u\|_{L^4}^2+\|\nabla d\cdot\Delta d\|_{L^4}^2\right)dt\nonumber\\
&\leq&C\int_0^T\left(\|\rho^{1/2}u_t\|_{L^2}^2+\|\nabla u_t\|_{L^2}^2+\|u\|_{H^2}^4+\|\nabla d\|_{H^2}^4\right)dt\nonumber\\
&\leq&C +C\int_0^T\|\nabla u_t\|_{L^2}^2dt\leq C.\label{3.37}
\end{eqnarray}
Moreover, by virtue Lemma \ref{lem2.1}, (\ref{3.28}) and (\ref{3.29}) we infer from (\ref{1.4}) that
\begin{eqnarray*}
\int_0^T\|\nabla^2d\|_{H^2}^2dt&\leq& C\int_0^T\left(\|d_t\|_{H^2}^2+\|u\cdot\nabla d\|_{H^2}^2+\||\nabla d|^2d\|_{H^2}^2\right)dt\\
&\leq&C\int_0^t\left(\|d_t\|_{H^2}^2+\|u\|_{H^2}^2\|\nabla d\|_{H^2}^2+\|\nabla d\|_{H^2}+\|\nabla d\|_{H^2}^6\right)dt\\
&\leq& C,
\end{eqnarray*}
where we have used the following Moser's type calculus inequality
(see \cite{Ma1984} ) that for $f,g\in H^s(\R^2)$ with $s\geq2$,
$$
\|fg\|_{H^s}\leq C\left(\|f\|_{L^\infty}\|g\|_{H^s}+\|f\|_{H^s}\|g\|_{L^\infty}\right)\leq C\|f\|_{H^s}\|g\|_{H^s}.
$$
This, together with (\ref{3.36}) and (\ref{3.37}), leads to (\ref{3.29}) immediately.\hfill$\square$

\vskip 2mm

The last step is to estimate the first and second order derivatives of the density.
\begin{lem}\label{lem3.5}
For every $0<T<T^*$, one has
\begin{eqnarray}
\sup_{0\leq t\leq T}\left(\|\nabla
\rho\|_{H^1}+\|\rho_t\|_{H^1}\right)+\int_0^T\|\nabla
u\|_{H^2}^2dt\leq C.\label{3.38}
\end{eqnarray}
\end{lem}
\pf
Differentiating (\ref{1.1}) with respect to $x_i$
$(i=1,2)$, multiplying  the resulting equation  by $|\nabla\rho|^{p-2}\partial_i
\rho$ with $p\geq 2$, and integrating it by parts
over $\mathbb{R}^2$, we obtain after summing up that
\begin{eqnarray*}
\frac{d}{ d t}\|\nabla\rho\|_{L^p}^p \leq C\|\nabla
u\|_{L^\infty}\|\nabla\rho\|_{L^p}^p\leq C\|\nabla
u\|_{W^{1,4}}\|\nabla\rho\|_{L^p}^p,
\end{eqnarray*}
which, combined with (\ref{3.29}) and Gronwall's inequality, yields
\begin{equation}
\|\nabla\rho\|_{L^p}^p\leq C\|\nabla\rho_0\|_{L^p}^p\exp\left(C\int_0^T\|\nabla u\|_{W^{1,4}}\dd t\right)\leq C,\quad\forall\; p\geq 2.\label{3.39}
\end{equation}

Similarly, by (\ref{3.38}) we also deduce from (\ref{1.1}) that
\begin{eqnarray*}
\frac{ d}{ d t}\|\nabla^2\rho\|_{L^2}^2&&\leq C\|\nabla
u\|_{L^\infty}\|\nabla^2\rho\|_{L^2}^2+C\|\nabla^2
u\|_{L^4}\|\nabla\rho\|_{L^4}\|\nabla^2\rho\|_{L^2}\nonumber\\
&&\leq C\|\nabla
u\|_{W^{1,4}}\left(1+\|\nabla^2\rho\|_{L^2}^2\right),
\end{eqnarray*}
so that
\begin{equation}
\|\nabla^2\rho\|_{L^2}^2\leq
C\left(1+\|\nabla^2\rho_0\|_{L^2}^2\right)\exp\left(C\int_0^T\|\nabla
u\|_{W^{1,4}}\dd t\right)\leq C.\label{3.40}
\end{equation}
As a result of (\ref{3.28}), (\ref{3.39}) and (\ref{3.40}), one
easily gets from (\ref{1.1}) that $\|\rho_t\|_{H^1}\leq C$.

Finally, it follows from (\ref{2.5}), (\ref{3.2}), (\ref{3.14}),
(\ref{3.28}) and (\ref{3.38}) that
\begin{eqnarray*}
\|\nabla u\|_{H^2}&\leq& C(\|\rho u_t\|_{H^1}+\|\rho u\cdot\nabla u\|_{H^1}+\|\nabla d\cdot\Delta d\|_{H^1})\nonumber\\
&\leq& C\left(\|\rho^{1/2} u_t\|_{L^2}+\|\nabla\rho\|_{L^4}\| u_t\|_{L^4}+\|\nabla u_t\|_{L^2}\right)\\
&&+C\left(\|u\|_{H^2}^2+\|u\|_{L^\infty}\|\nabla \rho\|_{L^4}\|\nabla u\|_{L^4}+\|\nabla d\|_{H^2}^2\right)\nonumber\\
&\leq& C\left(1+\|\nabla u_t\|_{L^2}\right),
\end{eqnarray*}
where we have also used Lemmas \ref{lem2.1} and \ref{lem2.2}.
Consequently,
$$
\int_0^T\|\nabla u\|_{H^2}^2dt\leq C+C\int_0^T\|\nabla u_t\|_{L^2}^2dt\leq C,
$$
which, together with (\ref{3.39}) and (\ref{3.40}), proves (\ref{3.38}).\hfill$\square$

\vskip 2mm

Collecting all the estimates in (\ref{3.2}), (\ref{3.3}) and Lemmas
\ref{lem3.1}--\ref{lem3.5} together, we arrive at (\ref{3.1}), and
hence, the proof of Theorem \ref{thm1.1} is complete.

\section{Proofs of Theorems \ref{thm1.2} and \ref{thm1.3}}
This section is concerned with the proofs of Theorems \ref{thm1.2}
and \ref{thm1.3}. We first prove Theorem \ref{thm1.2} by using
contradiction arguments. So, to do this, we assume otherwise that
\begin{equation}
\lim\limits_{T\to T^*}\int_0^{T}\|\nabla d\|_{L^r}^s dt\leq
M_0<\infty\label{4.1}
\end{equation}
with any $(r,s)$ satisfying (\ref{1.17}).

We begin the proof with the observation from the proof of Theorem
\ref{thm1.1} that, to remove the smallness condition (\ref{1.14})
and to obtain a global strong solution with generally large initial
data, it suffices to achieve the estimate of
$\|\nabla^2d\|_{L^2(0,T;L^2)}$ for any $0<T<T^*$. Moreover, it
follows from (\ref{3.9}) and (\ref{3.4}) that for any $0<T<T^*$,
\begin{eqnarray}
\int_0^T\|\nabla^2d\|_{L^2}^2dt&=&\int_0^T\|\Delta d+|\nabla
d|^2d\|_{L^2}^2dt+\int_0^T\|\nabla d\|_{L^4}^4dt\nonumber\\
&\leq& \frac{1}{2}E_0+\int_0^T\|\nabla d\|_{L^4}^4dt.\label{4.2}
\end{eqnarray}
Therefore, to bound $\|\nabla^2d\|_{L^2(0,T;L^2)}$, we only need to
deal with $\|\nabla d\|_{L^4(0,T;L^4)}$. This  will be done in the
following.

On one hand, assume that $(r,s)$ satisfies
\begin{equation}
\frac{1}{r}+\frac{1}{s}\leq \frac{1}{2}\quad{\rm with }\quad 4\leq
r\leq\infty.\label{4.3}
\end{equation}
Then, using H\"{o}lder inequality, (\ref{2.1}) and (\ref{3.4}), we
find
\begin{eqnarray}
\|\nabla d\|_{L^4}^4&=& \|\nabla d\|_{L^4}^2\|\nabla
d\|_{L^4}^2\nonumber\\
&\leq& C\|\nabla d\|_{L^2}^{(r-4)/(r-2)}\|\nabla d\|_{L^r}^{r/(r-2)} \|\nabla d\|_{L^2}\|\nabla^2d\|_{L^2}\nonumber\\
&\leq& C \|\nabla d\|_{L^r}^{r/(r-2)} \|\nabla^2d\|_{L^2}\nonumber\\
& \leq& \frac{1}{2}\|\nabla^2d\|_{L^2}^2+C\left(1+\|\nabla
d\|_{L^r}^s\right).\label{4.4}
\end{eqnarray}
Thus, putting (\ref{4.4}) into (\ref{4.2}) and using(\ref{4.1}), we
obtain
\begin{equation}
\int_0^T\|\nabla^2d\|_{L^2}^2dt \leq E_0+C\int_0^T\left(1+\|\nabla
d\|_{L^r}^s\right)dt\leq C,\label{4.5}
\end{equation}
provided $(r,s)$ satisfies (\ref{4.3}).

On the other hand, assume that $(r,s)$ satisfies
\begin{equation}
\frac{1}{r}+\frac{1}{s}\leq \frac{1}{2}\quad{\rm with }\quad 2<
r<4.\label{4.6}
\end{equation}
Then, by virtue of (\ref{2.2}), (\ref{3.4}) and H\"{o}lder inequality we find that
\begin{eqnarray}
\|\nabla d\|_{L^4}&\leq& C\|\nabla d\|_{L^r}^{1-\alpha}\|\nabla d\|_{L^p}^\alpha\nonumber\\
&\leq& C\|\nabla d\|_{L^r}^{1-\alpha}\|\nabla d\|^{2\alpha/p}_{L^2}\|\nabla^2 d\|^{\alpha(p-2)/p}_{L^2}\nonumber\\
&\leq& C\|\nabla d\|_{L^r}^{1-\alpha}\|\nabla^2
d\|^{\alpha(p-2)/p}_{L^2}\label{4.7}
\end{eqnarray}
with
$$
\frac{2r}{r-2}<p<\infty,\quad\alpha=\frac{(4-r)p}{4(p-r)}\in\left(0,\frac{1}{2}\right).
$$
As a result of (\ref{4.7}), we have by Young's inequality that
\begin{eqnarray}
\|\nabla d\|_{L^4}^4&\leq&C\|\nabla d\|_{L^r}^{4(1-\alpha)} \|\nabla^2 d\|_{L^2}^{4\alpha(p-2)/p}  \nonumber\\
&\leq& \frac{1}{2}\|\nabla^2d\|_{L^2}^2+C \|\nabla
d\|_{L^r}^{4(1-\alpha)p/(p-2\alpha(p-2))}\nonumber\\
&\leq& \frac{1}{2}\|\nabla^2d\|_{L^2}^2+C \left(1+\|\nabla
d\|_{L^r}^{s}\right),\label{4.8}
\end{eqnarray}
since direct calculations give
$$
\frac{4(1-\alpha)p}{p-2\alpha(p-2)}=\frac{2r}{r-2}\leq s.
$$
Thus, putting (\ref{4.8}) into (\ref{4.2}), by (\ref{4.1}) we also
obtain (\ref{4.5}), provided $(r,s)$ satisfies (\ref{4.6}).

To conclude, we have proved that there exists a positive constant $C$, depending on the initial data, $T^*$ and $M_0$, such that for any $0<T<T^*$,
$$
\int_0^T\|\nabla^2d\|_{L^2}^2 dt\leq C,
$$
provided (\ref{4.1}) holds. With the help of this and
(\ref{3.2})--(\ref{3.4}), following the arguments in the proofs of
Lemmas \ref{lem3.3}--\ref{lem3.5}, we arrive at (\ref{3.1}), which,
combined with the local existence theorem (see Lemma \ref{lem2.5}),
implies the solutions can be extended beyond $T^*$. This immediately
leads to a contradiction of $T^*$, and hence, the proof of Theorem
\ref{thm1.2} is complete.\hfill$\square$

\vskip 2mm

{\it Proof of Theorem \ref{thm1.2}.} In fact, by applying the
maximum principle to the equation of $d_3$ (i.e. the third component
of $d$), we have
$$
\inf_{x\in\R^2}d_3(x,t)\geq\inf_{x\in\R^2}d_{03}\geq\varepsilon,\quad\forall\;
t>0.
$$
So, it follows from (\ref{3.4}) and
(\ref{1.19}) that for any $0<T<\infty$,
\begin{equation}
\int_0^T\left(\|\Delta d\|_{L^2}^2+\|\nabla d\|_{L^4}^4\right)dt \leq
C(E_0,\varepsilon)\int_0^T\|\Delta d+|\nabla d|^2d\|_{L^2}^2dt\leq C(E_0,\varepsilon).
\end{equation}
This, combined with Theorem \ref{thm1.2} with $r=s=4$, implies that the strong solutions of  (\ref{1.1})--(\ref{1.4}), (\ref{1.10}) and
(\ref{1.11}) exist for all $T>0$. The proof of Theorem \ref{thm1.3} is thus finished.\hfill$\square$

\small{
}
\end{document}